\documentclass[bezier]{article}
\usepackage{graphicx}
\usepackage{amsmath,amssymb,amsfonts, euscript}
\newtheorem{theorem}{Theorem}[section]

\newtheorem{lemma}[theorem]{Lemma}

\newtheorem{example}[theorem]{Example}

\title{\bf\Large A class of multiplicative splitting iterations for solving the continuous Sylvester equation}

\author{Yu Huang$^1$, Mohammad Khorsand Zak$^2$\footnote{Corresponding author, {\it email}: mo.khorsand@mail.um.ac.ir} and Emran Tohidi$^3$
\\{\footnotesize{\it $^1$College of Mathematics and Statistics, Nanjing University of Information Science and Technology, Nanjing 210044, P. R. China}}
\\{\footnotesize{\it $^2$Department of Applied Mathematics, Aligudarz Branch, Islamic Azad University, Aligudarz, Iran}}
\\{\footnotesize{\it $^3$Department of Mathematics, Kosar University of Bojnord, P.O. Box 9415615458, Bojnord, Iran}}
}

\date{}

\begin{document}
\maketitle

\begin{abstract}
For solving the continuous Sylvester equation, a class of the multiplicative splitting iteration method is presented. We consider two symmetric positive definite splittings for each coefficient matrix of the continuous Sylvester equations and it can be equivalently written as two multiplicative splitting matrix equations.
When both coefficient matrices in the continuous Sylvester equation are (non-symmetric) positive semi-definite,
and at least one of them is positive definite; we can choose Hermitian and skew-Hermitian (HS) splittings of matrices $A$ and $B$, in the first equation, and the splitting of the Jacobi iterations for matrices $A$ and $B$, in the second equation in the multiplicative splitting iteration method.
Convergence conditions of this method are studied and numerical experiments show the efficiency of this method.

\end{abstract}

\section{Introduction}
Matrix equations arise in a number of problems of scientific computations and engineering applications, such as control theory \cite{Datta}, model reduction \cite{Baur}, signal and image processing \cite{Bouhamidi} and many researchers focus on the matrix equations \cite{Bai2, Beik, Dehghan, Hajarian, Jbilou, Simoncini}. Nowadays, the continuous Sylvester equation is possibly the most famous and the most broadly employed linear matrix equation \cite{Bai2, Bouhamidi, Dehghan1, El, El1, Hajarian1, Khorsand2, Khorsand3, Simoncini, Tohidi}, and is given as
\begin{equation}\label{1}
AX+XB=C,
\end{equation}
where $A\in\mathbb{R}^{n\times n}$, $B\in\mathbb{R}^{m\times m}$ and
$C\in\mathbb{R}^{n\times m}$ are defined matrices and $X\in\mathbb{R}^{n\times m}$ is an unknown matrix.
The continuous Sylvester equation \eqref{1} has a unique solution if and only if $A$ and $-B$ have no common eigenvalues, which will be assumed throughout this paper.

In general, the dimensions of $A$ and $B$ may be orders of magnitude different, and this fact is key in selecting the most appropriate numerical solution strategy \cite{Simoncini}.
For solving general Sylvester equations of small size we use some methods which classified such as direct methods. Some of these direct methods are the Bartels-Stewart \cite{Bartels} and the Hessenberg-Schur \cite{Golub} methods which consist of transforming coefficient matrices $A$ and $B$ into triangular or Hessenberg form by an orthogonal similarity transformation and then solving the resulting system directly by a back-substitution process.
When the coefficient matrices $A$ and $B$ are large and sparse, iterative methods are often the methods of choice for solving the Sylvester equation \eqref{1} efficiently and accurately. Many iterative methods were developed for solving matrix equations, such as the alternating direction implicit (ADI) method \cite{Benner}, the Krylov subspace based algorithms \cite{Hu, Khojasteh, El}, the Hermitian and skew-Hermitian splitting (HSS) method, and the inexact variant of HSS (IHSS) iteration method \cite{Bai2}, The nested splitting conjugate gradient (NSCG) method \cite{Khorsand1, Khorsand2} and the nested splitting CGNR (NS-CGNR) method \cite{Khorsand3}.

In order to study the numerical methods, we often rewrite the continuous Sylvester equation \eqref{1} as the linear system of equations
\begin{equation}\label{2}
\mathcal{A}x=c,
\end{equation}
where the matrix $\mathcal{A}$ is of dimension $nm\times nm$ and is
given by
\begin{equation}\label{3}
\mathcal{A}=I_m\otimes A+B^T\otimes I_n,
\end{equation}
where $\otimes$ denotes the Kronecker product $(A\otimes B=[a_{ij}B])$ and
$$
\begin{array}{l}
c=vec(C)=(c_{11},c_{21},\cdots,c_{n1},c_{12},c_{22},\cdots,c_{n2},\cdots,c_{nm})^T\\
x=vec(X)=(x_{11},x_{21},\cdots,x_{n1},x_{12},x_{22},\cdots,x_{n2},\cdots,x_{nm})^T.
\end{array}
$$

Of course, this is quite expensive and a numerically poor way to determine the solution $X$ of the continuous Sylvester equation \eqref{1}, as the linear system of equations \eqref{2} is costly to solve and can be ill-conditioned.

Now, we recall some necessary notations and useful results, which will be used in the following section.
In this paper, we use $\lambda(M), ||M||_2$, $||M||_F$ and $I_n$ to denote the eigenvalue, the spectral norm, the Frobenius norm of a matrix $M\in\mathbb{R}^{n\times n}$, and the identity matrix with dimension $n$, respectively. Note that $||.||_2$ is also used to represent the 2-norm of a vector.
For nonsingular matrix $\mathcal{B}$, we denote by $\kappa(\mathcal{B})=||\mathcal{B}||_2||\mathcal{B}^{-1}||_2$ its spectral condition number, and for a symmetric a positive definite matrix $\mathcal{B}$, we define the $||\cdot||_{\mathcal{B}}$ norm of a vector $x\in\mathbb{R}^n$ as $||x||_{\mathcal{B}}=\sqrt{x^H\mathcal{B}x}$. Then the induced $||\cdot||_{\mathcal{B}}$ norm of a matrix $M\in\mathbb{R}^{n\times n}$ is define as $||M||_{\mathcal{B}}=||\mathcal{B}^{\frac{1}{2}}M\mathcal{B}^{-\frac{1}{2}}||_2$. In addition it holds that $||Mx||_{\mathcal{B}}\leq||M||_{\mathcal{B}}||x||_{\mathcal{B}}$, $||M||_{\mathcal{B}}\leq\sqrt{\kappa(\mathcal{B})}||M||_2$ and $||I||_{\mathcal{B}}=1$, where $I$ is the identity matrix.
For any matrices $A=[a_{ij}]$ and $B=[b_{ij}]$, $A\otimes B$ denotes the Kronecker product defined as $A\otimes B=[a_{ij}B]$. For the matrix $X=(x_1, x_2,\cdots,x_m)\in\mathbb{R}^{n\times m}$, $vec(X)$ denotes the $vec$ operator defined as $vec(X)=(x_1^T,x_2^T,\cdots,x_m^T)^T$.
Moreover, for a matrix $M\in\mathbb{R}^{n\times n}$ and the vector $vec(M)\in\mathbb{R}^{nm}$, we have $||M||_F=||vec(M)||_2$.

For matrix $\mathcal{A}\in\mathbb{R}^{n\times n}$, $\mathcal{A}=\mathcal{B}-\mathcal{C}$ is called a splitting of the matrix $\mathcal{A}$ if $\mathcal{B}$ is nonsingular. This splitting is a convergent splitting if $\rho(\mathcal{B}^{-1}\mathcal{C})<1$; and a contractive splitting if $||\mathcal{B}^{-1}\mathcal{C}||<1$ for some matrix norm.

The reminder of this paper is organized as follows. Section \ref{main} presents our main contribution. In other words, the multiplicative splitting iteration (MSI) method for the continuous Sylvester equation and its convergence properties are studied deeply. Section \ref{numer} is devoted to an extensive numerical experiments with full comparison with other state of the art methods in the literature. Finally, we present our conclusions in Section \ref{conc}.

\section{Multiplicative splitting iterations}\label{main}
\subsection{Traditional MSI method}
Consider the linear system of equations \eqref{2}. Let $\mathcal{A}=\mathcal{M}_i-\mathcal{N}_i$ $(i=1,2)$ be two splittings of the coefficient matrix $\mathcal{A}$. The MSI method for solving the system of linear equations \eqref{2} is defined as follows \cite{Bai1}:\\
\textbf{MSI method for linear system of equations}\\
\verb"Given an initial guess" $x^{(0)}\in\mathbb{R}^{n}$,\\
\verb"For" $k=1,2,\cdots$ \verb"until convergence, do"
\begin{quote}
    $u^{(k+1)}=\mathcal{M}^{-1}_1\mathcal{N}_1x^{(k)}+\mathcal{M}^{-1}_1c$\\
    $x^{(k+1)}=\mathcal{M}^{-1}_2\mathcal{N}_2u^{(k+1)}\mathcal{M}^{-1}_2c$
\end{quote}
\verb"end"

The MSI method can be equivalently written in the form
\[
x^{(k+1)}=\mathcal{T}_{msi}x^{(k)}+\mathcal{G}_{msi}c, \hspace{1cm} k=0,1,2,\cdots
\]
where $\mathcal{T}_{msi}=\mathcal{M}^{-1}_2\mathcal{N}_2\mathcal{M}^{-1}_1\mathcal{N}_1$ and
$\mathcal{G}_{msi}=\mathcal{M}^{-1}_2\mathcal{N}_2\mathcal{M}^{-1}_1+\mathcal{M}^{-1}_2$. See \cite{Bai1} for more details.
\subsection{MSI method for the Sylvester equation}
Based on the MSI method proposed in \cite{Bai1}, we obtain the MSI method for the continuous Sylvester equation.
Let $A=M_i-N_i$ and $B=P_i-Q_i,~ (i=1,2)$ be two splittings of the matrices $A$ and $B$, such that $M_i$ and $P_i,~(i=1,2)$ are symmetric positive definite. The continuous Sylvester equation (\ref{1}) can be equivalently written as the multiplicative splitting matrix equations
\[
\left\{
\begin{array}{l}
  M_1U+UP_1=N_1X+XQ_1+C \\
  M_2X+XP_2=N_2U+UQ_2+C
\end{array}
\right.
\]
Under the assumption that $M_i$ and $P_i,~(i=1,2)$ are symmetric positive definite, we easily know that there is no common eigenvalues between the matrices $M_i$ and $-P_i,~(i=1,2)$, so that this two multiplicative splitting matrix equations have unique solutions for all given right hand side matrices.

Now, based on the above observations, we can establish the following multiplicative splitting iterations for solving the continuous Sylvester equation \eqref{1}:\\
\textbf{MSI method for Sylvester equation }\\
\verb"Given an initial guess" $X^{(0)}\in\mathbb{R}^{m\times n}$,\\
\verb"For" $k=1,2,\cdots$ \verb"until convergence, do"
\begin{quote}
    \verb"Solve" $M_1U^{(k+1)}+U^{(k+1)}P_1=N_1X^{(k)}+X^{(k)}Q_1+C$\\
    \verb"Solve" $M_2X^{(k+1)}+X^{(k+1)}P_2=N_2U^{(k+1)}+U^{(k+1)}Q_2+C$
\end{quote}
\verb"end"

In special case, when both coefficient matrices $A$ and $B$, in Sylvester equation \eqref{1} are (non-symmetric) positive semi-definite,
and at least one of them is positive definite; we can choose Hermitian and skew-Hermitian (HS) splittings of matrices $A$ and $B$, in the first equation in MSI method, and the splitting of the Jacobi iterations \cite{Saad} for matrices $A$ and $B$, in the second equation in MSI method. Therefore, we can rewrite this method as following: \\
\verb"Given an initial guess" $X^{(0)}\in\mathbb{R}^{m\times n}$,\\
\verb"For" $k=1,2.\cdots$ \verb"until convergence, do"
\begin{quote}
    \verb"Solve system" $H_AU^{(k+1)}+U^{(k+1)}H_B=S_AX^{(k)}+X^{(k)}S_B+C$\\
    \verb"Solve system" $D_AX^{(k+1)}+X^{(k+1)}D_B=N_AU^{(k+1)}+U^{(k+1)}N_B+C$
\end{quote}
\verb"end"

Achieving to two Sylvester equations that we can easily solve them, is our motivation for choice of these splittings. Because the first system can be solved by Sylvester conjugate gradient method \cite{Evans}, and the following routine can be used for direct solution of the second system.\\
\textbf{Directly solution of matrix equation $D_AX+XD_B=C$ }\\
\verb"For" $i=1:n$
\begin{quote}
\verb"For" $j=1:m$
\begin{quote}
$x_{ij}=\frac{c_{ij}}{a_{ii}+b_{jj}}$
\end{quote}
\verb"end"
\end{quote}
\verb"end"

\subsection{Convergence analysis}
In the sequel, we need the following lemmas.
\begin{lemma}\label{l1}\cite{Bai1}
Let $B, C\in\mathbb{R}^{n\times n}$ be two Hermitian matrices. Then $BC=CB$ if and only if $B$ and $C$ have a common set of orthonormal eigenvectors.
\end{lemma}
\begin{lemma}\label{lema1}\cite{Kelley}
Let $\mathcal{A}\in\mathbb{R}^{n\times n}$ be a symmetric positive definite matrix. Then for all $x\in\mathbb{R}^n$, we have $||\mathcal{A}^{\frac{1}{2}}x||_2=||x||_\mathcal{A}$ and
$$\sqrt{\lambda_{\min}(\mathcal{A})}||x||_\mathcal{A}\leq||\mathcal{A}x||_2\leq\sqrt{\lambda_{\max}(\mathcal{A})}||x||_\mathcal{A}.$$
\end{lemma}
\begin{lemma}\label{lema2}\cite{Lutkepohl}
Suppose that $A,B\in\mathbb{R}^{n\times n}$ be two Hermitian matrices, and denote
the minimum and the maximum eigenvalues of a matrix $M$ with
$\lambda_{\min}(M)$ and $\lambda_{\max}(M)$, respectively. Then
\[
\begin{array}{l}
\lambda_{\max}(A+B)\leq\lambda_{\max}(A)+\lambda_{\max}(B),\\
\lambda_{\min}(A+B)\geq\lambda_{\min}(A)+\lambda_{\min}(B).
\end{array}
\]
\end{lemma}
\begin{lemma}\label{lema3}\cite{Lutkepohl}
Let $A,B\in\mathbb{C}^{n\times n}$, and $\lambda$ and $\mu$ be the eigenvalues of $A$ and $B$, and $x$ and $y$ be the corresponding
eigenvectors, respectively. Then $\lambda\mu$ is an eigenvalue of $A\otimes B$ corresponding to the eigenvector $x\otimes y$.
\end{lemma}
\begin{lemma}\label{lema4}
Suppose that $\mathcal{A}=\mathcal{M}-\mathcal{N}$ is a splitting such that $\mathcal{M}$ is symmetric positive definite, with $\mathcal{M}=I_m\otimes M+P^T\otimes I_n$ and $\mathcal{N}=I_m\otimes N+Q^T\otimes I_n$.
If
\[
\theta^3\frac{\max|\lambda(N)|+\max|\lambda(Q)|}{\lambda_{\min}(M)+\lambda_{\min}(P)}<1,
\]
where $\theta=\sqrt{\frac{\lambda_{\max}(M)+\lambda_{\max}(P)}{\lambda_{\min}(M)+\lambda_{\min}(P)}}$,
then $||\mathcal{M}^{-1}\mathcal{N}||_\mathcal{M}<1$.
\end{lemma}
\textrm{proof.} By Lemmas \ref{lema2} and \ref{lema3}, we have
\[
||\mathcal{M}||_2=\lambda_{\max}(\mathcal{M})\geq\lambda_{\min}(\mathcal{M})\geq\lambda_{\min}(M)+\lambda_{\min}(P),
\]
and
\[
||\mathcal{N}||_2=\max_{\lambda\in\Lambda(\mathcal{N})}|\lambda(\mathcal{N})|\leq\max|\lambda(N)|+\max|\lambda(Q)|,
\]
Therefore, it follows that
\[
\begin{array}{rcl}
||\mathcal{M}^{-1}\mathcal{N}||_{\mathcal{M}}& \leq &\sqrt{\kappa(\mathcal{M})}||\mathcal{M}^{-1}\mathcal{N}||_2\\
                 & \leq& \sqrt{\kappa(\mathcal{M})}||\mathcal{M}^{-1}||_2||\mathcal{N}||_2\\
                 & \leq & (\kappa(\mathcal{M}))^{\frac{3}{2}}\frac{||\mathcal{N}||_2}{||\mathcal{M}||_2}\\
                 & \leq & (\kappa(\mathcal{M}))^{\frac{3}{2}}\frac{\max|\lambda(N)|+\max|\lambda(Q)|}{\lambda_{\min}(M)+\lambda_{\min}(P)}.
\end{array}
\]
Again, the use of Lemmas \ref{lema2} and \ref{lema3} implies that
\begin{equation}\label{8a}
\sqrt{\kappa(\mathcal{M})}=
\sqrt{\frac{\lambda_{\max}(\mathcal{M})}{\lambda_{\min}(\mathcal{M})}}\leq
\sqrt{\frac{\lambda_{\max}(M)+\lambda_{\max}(P)}{\lambda_{\min}(M)+\lambda_{\min}(P)}}=\theta.
\end{equation}
So, we can write
\begin{equation}\label{8b}
||\mathcal{M}^{-1}\mathcal{N}||_{\mathcal{M}}\leq\theta^3\frac{\max|\lambda(N)|+\max|\lambda(Q)|}{\lambda_{\min}(M)+\lambda_{\min}(P)}.
\end{equation}
This clearly proves the lemma.
\begin{theorem}
Let $A\in\mathbb{R}^{n\times n}$ and $B\in\mathbb{R}^{m\times m}$ and consider two splittings  $A=M_i-N_i$ and $B=P_i-Q_i~(i=1,2)$ such that $M_i$ and $P_i,~(i=1,2)$ are symmetric positive definite. Denote by $\mathcal{A}=\mathcal{M}_i-\mathcal{N}_i~(i=1,2)$ with $\mathcal{M}_i=I_m\otimes M_i+P^T_i\otimes I_n$ and $\mathcal{N}_i=I_m\otimes N_i+Q^T_i\otimes I_n~(i=1,2)$, and assume that $\mathcal{M}_1\mathcal{A}^{-1}$ and $\mathcal{M}_2\mathcal{A}^{-1}$ are Hermitian matrices and $\mathcal{M}_1\mathcal{A}^{-1}\mathcal{M}_2=\mathcal{M}_2\mathcal{A}^{-1}\mathcal{M}_1$. Then the MSI method is convergent if $\varrho_1\varrho_2<1$, where
\[
\varrho_i=\theta_i^3\frac{\max|\lambda(N_i)|+\max|\lambda(Q_i)|}{\lambda_{\min}(M_i)+\lambda_{\min}(P_i)},~and~
\theta_i=\sqrt{\frac{\lambda_{\max}(M_i)+\lambda_{\max}(P_i)}{\lambda_{\min}(M_i)+\lambda_{\min}(P_i)}},
~(i=1,2).
\]
\end{theorem}
\textrm{Proof.} By making use of the Kronecker product, we can rewrite the above described MSI method in the following matrix-vector form:
\[
\left\{
\begin{array}{l}
  (I_m\otimes M_1+P^T_1\otimes I_n)u^{(k+1)}=(I_m\otimes N_1+Q^T_1\otimes I_n)x^{(k)}+c \\
  (I_m\otimes M_2+P^T_2\otimes I_n)x^{(k+1)}=(I_m\otimes N_2+Q^T_2\otimes I_n)u^{(k+1)}+c
\end{array}
\right.
\]
which can be arranged equivalently as
\[
\left\{
\begin{array}{l}
  \mathcal{M}_1u^{(k+1)}=\mathcal{N}_1x^{(k)}+c \\
  \mathcal{M}_2x^{(k+1)}=\mathcal{N}_2u^{(k+1)}+c
\end{array}
\right.
\]
which can be obtained the following iteration method

\begin{equation}\label{MSI}
\left\{
\begin{array}{l}
  u^{(k+1)}=\mathcal{M}^{-1}_1\mathcal{N}_1x^{(k)}+\mathcal{M}^{-1}_1c \\
  x^{(k+1)}=\mathcal{M}^{-1}_2\mathcal{N}_2u^{(k+1)}+\mathcal{M}^{-1}_2c
\end{array}
\right.
\end{equation}

Evidently, the above iteration scheme is the MSI-method \cite{Bai1} for solving system of linear equations \eqref{2} with $\mathcal{A}=\mathcal{M}_i-\mathcal{N}_i~(i=1,2)$. The MSI iteration \eqref{MSI} can be neatly expressed as a stationary fixed-point iteration as follows,
$$x^{(k+1)}=\mathcal{T}x^{(k)}+\mathcal{G}c$$
with $\mathcal{T}=\mathcal{M}^{-1}_2\mathcal{N}_2\mathcal{M}^{-1}_1\mathcal{N}_2$ and $\mathcal{G}=\mathcal{M}^{-1}_2\mathcal{N}_2\mathcal{M}^{-1}_1+\mathcal{M}^{-1}_1$.

Because $\mathcal{M}_1\mathcal{A}^{-1}\mathcal{M}_2=\mathcal{M}_2\mathcal{A}^{-1}\mathcal{M}_1$ is equivalent to that the two matrices $\mathcal{M}_1\mathcal{A}^{-1}$ and $\mathcal{M}_2\mathcal{A}^{-1}$ are commutative, according to Lemma \ref{l1} we know that $\mathcal{M}_1\mathcal{A}^{-1}$ and $\mathcal{M}_2\mathcal{A}^{-1}$ have a common set of orthonormal eigenvectors. That is say, there exists a unitary matrix $\mathcal{Q}\in\mathbb{R}^{nm\times nm}$ and two diagonal matrices $\Lambda_i=\textit{diag}(\lambda_1^{(i)},\lambda_2^{(i)},\cdots,\lambda_{nm}^{(i)}),~i=1,2$, such that $\mathcal{Q}\mathcal{M}_i^{-1}\mathcal{A}\mathcal{Q}^{\ast}=\Lambda_i,~i=1,2$. Noticing that
\[
\begin{array}{rl}
\mathcal{T}&=\mathcal{M}_2^{-1}\mathcal{N}_2\mathcal{M}_1^{-1}\mathcal{N}_1\\
 &=\mathcal{M}_2^{-1}(\mathcal{M}_2-\mathcal{A})\mathcal{M}_1^{-1}(\mathcal{M}_1-\mathcal{A})\\
 &=(I-\mathcal{M}_2^{-1}\mathcal{A})(I-\mathcal{M}_1^{-1}\mathcal{A})\\
 &=(\mathcal{Q}^{\ast}\mathcal{Q}-\mathcal{Q}^{\ast}\mathcal{Q}\mathcal{M}_2^{-1}\mathcal{A}\mathcal{Q}^{\ast}\mathcal{Q})(\mathcal{Q}^{\ast}\mathcal{Q}-\mathcal{Q}^{\ast}\mathcal{Q}\mathcal{M}_1\mathcal{A}\mathcal{Q}^{\ast}Q)\\
 &=(\mathcal{Q}^{\ast}\mathcal{Q}-\mathcal{Q}^{\ast}\Lambda_2\mathcal{Q})(\mathcal{Q}^{\ast}\mathcal{Q}-\mathcal{Q}^{\ast}\Lambda_1\mathcal{Q})\\
 &=\mathcal{Q}^{\ast}(I-\Lambda_2)\mathcal{Q}^{\ast}\mathcal{Q}(I-\Lambda_1)\mathcal{Q}\\
 &=\mathcal{Q}^{\ast}(I-\Lambda_2)(I-\Lambda_1)\mathcal{Q}
\end{array}
\]
so we have
\[
\begin{array}{rl}
\rho(\mathcal{T})&\leq\max_{1\leq i,j\leq nm}|(1-\lambda_i^{(2)})(1-\lambda_j^{(1)})|\\
       &\leq\max_{1\leq i\leq nm}|(1-\lambda_i^{(2)})|\max_{1\leq j\leq nm}|(1-\lambda_j^{(1)})|\\
       &=\rho(I-\mathcal{M}_2^{-1}\mathcal{A})\rho(I-\mathcal{M}_1^{-1}\mathcal{A})\\
       &=\rho(\mathcal{M}_2^{-1}\mathcal{N}_2)\rho(\mathcal{M}_1^{-1}\mathcal{N}_1)\\
       &\leq||\mathcal{M}_2^{-1}\mathcal{N}_2||_{\mathcal{M}_2}||\mathcal{M}_1^{-1}\mathcal{N}_1||_{\mathcal{M}_1}
\end{array}
\]

Therefore, by Lemma \ref{lema4} we have
\[
\rho(\mathcal{T})\leq
\theta_1^3\frac{\max|\lambda(N_1)|+\max|\lambda(Q_1)|}{\lambda_{\min}(M_1)+\lambda_{\min}(P_1)}
\theta_2^3\frac{\max|\lambda(N_2)|+\max|\lambda(Q_2)|}{\lambda_{\min}(M_2)+\lambda_{\min}(P_2)}
=\varrho_1\varrho_2
\]
and this completes the proof.

\section{Numerical results}\label{numer}
All numerical experiments presented in this section were computed in double precision with a number of MATLAB codes. All
iterations are started from the zero matrix for initial $X^{(0)}$ and terminated when the current iterate satisfies
$\frac{\|R^{(k)}\|_F}{\|R^{(0)}\|_F}\leq10^{-8}$, where $R^{(k)}=C-AX^{(k)}-X^{(k)}B$ is the residual of the $k$th iterate. Also we use the tolerance $\varepsilon=0.01$ for inner iterations in corresponding methods.
For each experiment we report the number of iterations or the number of total outer iteration steps and CPU time, and compare the MSI method with NSCG \cite{Khorsand2}, GMRES \cite{Khojasteh}, BiCGSTAB \cite{El1} and HSS \cite{Bai2} iterative methods.
\begin{example}\label{Ex1}
 For this example, we use the matrices
\[
A=B=M+2rN+\frac{100}{(n+1)^2}I,
\]
where $M={\rm tridiag}(-1,2,-1)$, $N={\rm tridiag}(0.5,0,-0.5)$ and $r=0.01$ \cite{Bai2}.
\end{example}
We apply the iteration methods to this problem with different dimensions. The results are given in Table \ref{Tab1}.

\begin{table}[h]
\centering
\caption{\label{Tab1}Results of the Example \ref{Ex1}}
{\begin{tabular}{@{}lrrrrr}\hline
Method &{\centering $(32,32)$}&{\centering $(64,64)$} &{\centering $(128,128)$}&{\centering $(256,256)$}&{\centering $(512,512)$}\\
\hline
MSI        & (4,60,0.04)  & (5,155,0.23)  & (6,385,1.76)   & (7,910,9.62)     & (11,3026,144.62) \\
NSCG       & (4,62,0.02)  & (5,152,0.08)  & (6,384,1.18)   & (7,899,7.89)     & (11,3025,123.96) \\
HSS        & (48,576,0.17)& (89,1513,0.98)&(164,3662,13.98)& (298,7464,89.1)  & (541,15892,1429.6)\\
GMRES      & (7,70,0.05)  & (17,170,0.56) & (52,520,5.59)  & (178,1780,49.7)  & (610,6100,1252.7)\\
BiCGSTAB   & (39,-,0.02)  & (74,-,0.15)   & (143,-,2.35)   & (277,-,19.53)    & (635,-,359.5) \\
\hline
\end{tabular}}
\end{table}

The pair $(n,m)$ in the first row of the Table \ref{Tab1}, represents the dimension of matrices $A$ and $B$, respectively.
Moreover, the triplex $(a, b, c)$ in Table \ref{Tab1}, represents the number of outer iterations, the number of total iterations and the CPU time (in seconds), respectively.
From the results presented in the Table \ref{Tab1}, it can be seen that for large dimensions, the MSI and the NSCG methods are more efficient than the other methods.

\begin{example}\label{Ex2}
We consider the continuous Sylvester equation \eqref{1} with $n=m$ and the coefficient matrices
\[
\left\{
\begin{array}{l}
  A=\text{diag}(1, 2, \cdots, n)+rL^T, \\
  B=2^{-t}I_n+\text{diag}(1, 2, \cdots, n)+rL^T+2^{-t}L,
\end{array}
\right.
\]
with $L$ the strictly lower triangular matrix having ones in the lower triangle part \cite{Bai2}.
\end{example}
The iteration methods were used for this problem and the results are given in Table \ref{Tab2}. Moreover, we compare the convergence history of the iterative methods by residual norm decreasing in Figure \ref{fig1}.
\begin{table}[h]
\centering
\caption{\label{Tab2}\small{Results of the Example \ref{Ex2} }}
{\begin{tabular}{@{}lcccc}\hline
Method   &{\centering out-itr}&{\centering CPU time} &{\centering res-norm } \\
\hline
MSI      & 5       & 16.43  & 0.0029 \\
NSCG     & 8       & 21.65  & 0.0070 \\
HSS      & 99      &326.71  & 0.0288 \\
GMRES(10)& 20      &49.87   & 0.0027 \\
BiCGSTAB & 75      & 16.22  & 0.0028 \\
\hline
\end{tabular}}
\end{table}
In Table \ref{Tab2}, we report the number of outer iterations (out-itr), the CPU time and the residual norm (res-norm) after convergence. For this example, we observe that the MSI method is superior to the other iterative methods in terms of the number of iterations and it is similar to the BiCGSTAB method in terms of CPU time. Comparing the convergence history of the iterative methods by residual norm decreasing shows that the MSI method converges more rapid and smooth than the BiCGSTAB method ( see Figure \ref{fig1} ).

\begin{figure}\label{fig1}\center
  \includegraphics[scale=0.4]{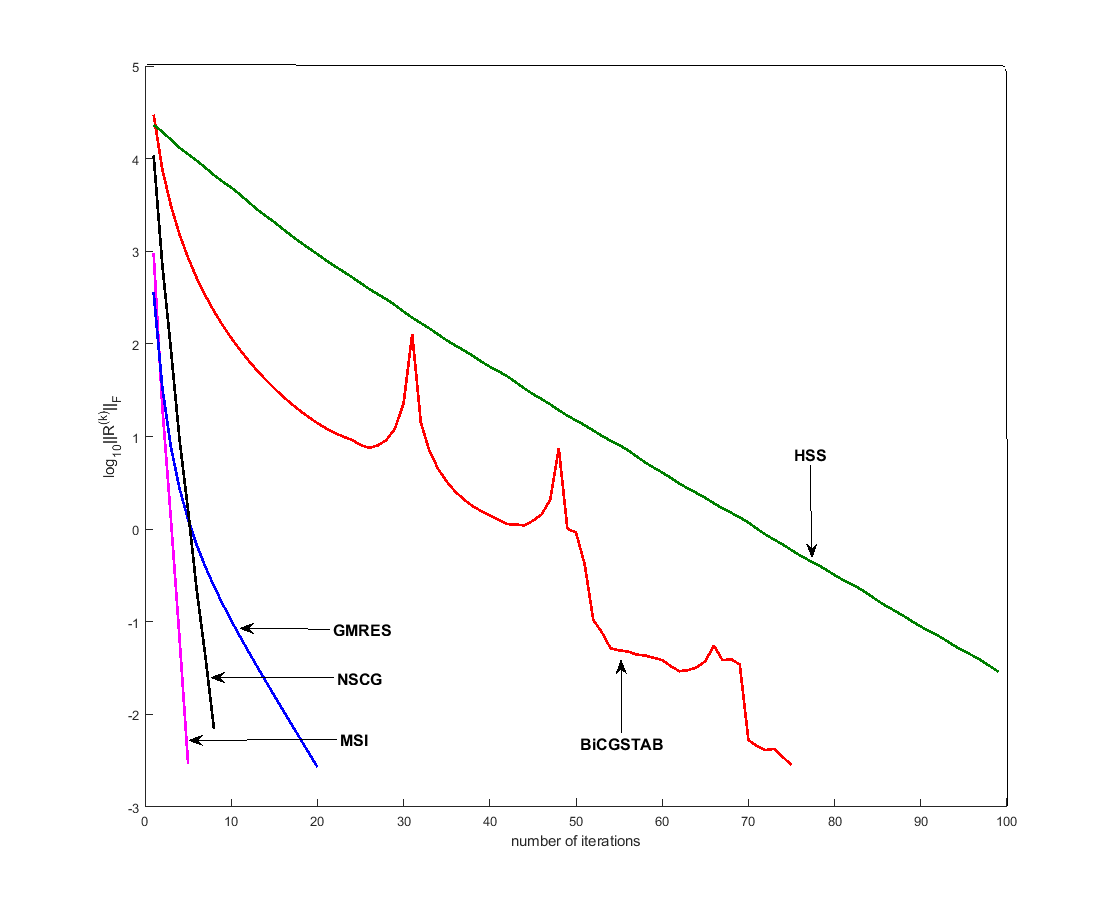}
  \caption{\label{fig1}{\small Convergence history of MSI versus the other iterative methods for Example \ref{Ex2} }}
\end{figure}

\begin{example}\label{Ex3}
For this example, we used the nonsymmetric sparse matrix SHERMAN3 of dimension $5005\times 5005$ with $20033$ nonzero entries from the Harwell-Boeing collection \cite{Duff} instead the coefficient matrix $A$. For the coefficient matrix $B$,  we used
$B={\rm tridiag}(-1,4,-2)$ of dimension $8\times8$ \cite{Khorsand2}.
\end{example}
We apply the iteration methods to this problem and the results are given in Table \ref{Tab3}. Moreover, we compare the convergence history of the iterative methods by residual norm decreasing in Figure \ref{fig2}.
\begin{table}[h]
\centering
\caption{\label{Tab3}\small{Results of the Example \ref{Ex3} }}
{\begin{tabular}{@{}lcccc}\hline
Method   &{\centering out-itr}&{\centering CPU time} &{\centering res-norm } \\
\hline
MSI      & 34      & 78.437  & 1.57e-4 \\
NSCG     & 64      & 121.265 & 2.61e-4 \\
HSS      & $>$5000 &$>$1000  & 2.32 \\
GMRES(10)& $>$5000 &$>$1000  & 247.77 \\
BiCGSTAB &  $\dag$ & $\dag$  &   NaN  \\
\hline
\end{tabular}}
\end{table}
In Table \ref{Tab3}, we report the number of outer iterations (out-itr), the CPU time and the residual norm (res-norm) after convergence or in 5000 outer iterations. For this example, we observe that the MSI method is superior to the other iterative methods in terms of the number of iterations and CPU times, the NSCG method has an acceptable performance. Furthermore, the HSS and the GMRES methods have a very slow convergence rate, and the BiCGSTAB method was diverged ( see Figure \ref{fig2} ).
\begin{figure}\label{fig2}\center
  \includegraphics[scale=0.4]{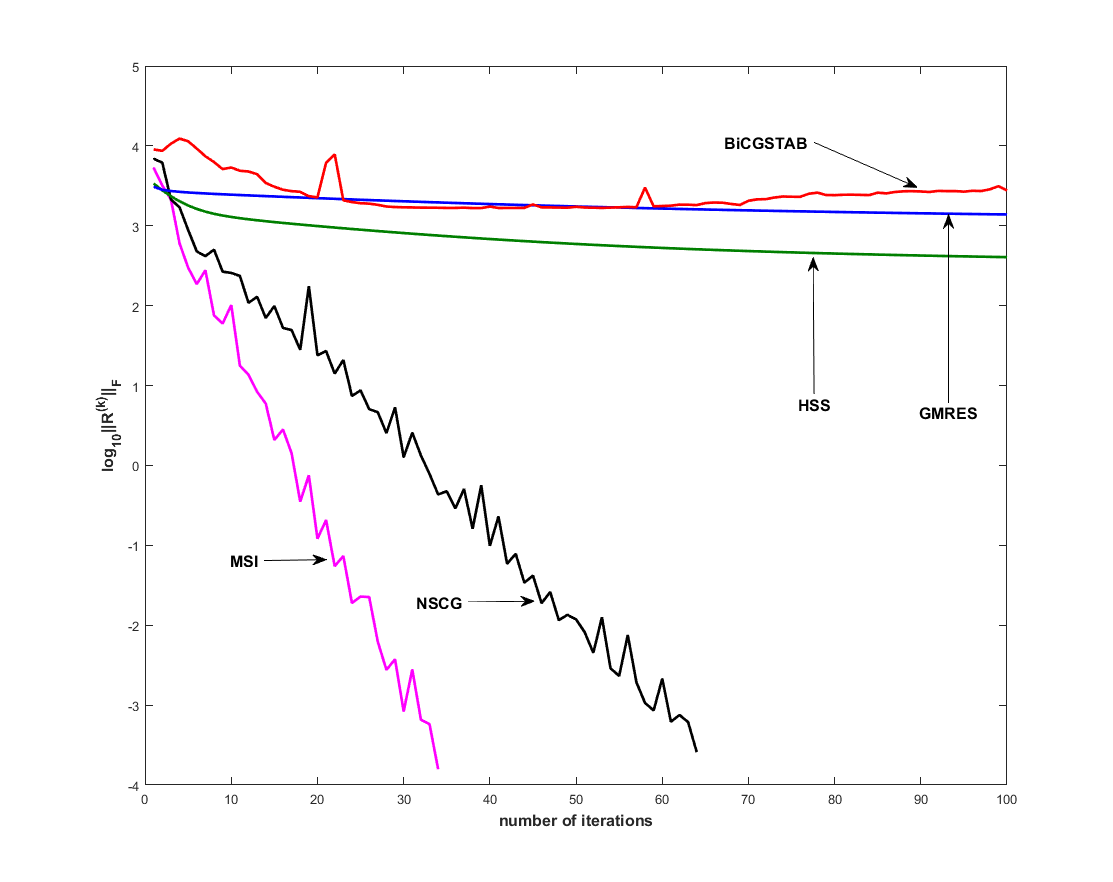}
  \caption{\label{fig2}{\small Convergence history of MSI versus the other iterative methods for Example \ref{Ex3} }}
\end{figure}

\section{Conclusion}\label{conc}
In this paper, we have proposed an efficient iterative method for solving the continuous Sylvester equation $AX+XB=C$. This method employs
two symmetric positive definite splittings of the coefficient matrices $A$ and $B$ and present a multiplicative splitting iteration method.

We have compared the MSI method with a well-known iterative methods such as the NSCG method, the HSS method, the BiCGSTAB method and the GMRES method for some problems. We have observed that, for these problems the MSI method is more efficient versus the other methods.

\textbf{Acknowledgments} Work of the first author Yu Huang was supported by National Natural Science Founding of China (No. 11771214 and No. 11801276).


\end{document}